\documentclass{amsart}

\usepackage[utf8]{inputenc}
\usepackage{amssymb,amsthm,amsmath}
\usepackage[foot]{amsaddr}
\usepackage{graphics,psfrag,graphicx,subfigure,epsfig}
\usepackage{enumerate}
\usepackage{hyperref}
\usepackage[all]{xy}
\usepackage{mathtools} 
\usepackage{tikz-cd}
\usepackage{amsrefs}
\usepackage{euscript}

\newtheorem{thm}{Theorem}[section]
\newtheorem{prop}[thm]{Proposition}

\newtheorem{theorem}[thm]{Theorem}
\newtheorem{lem}[thm]{Lemma}
\newtheorem{proposition}[thm]{Proposition}

\newtheorem{cor}[thm]{Corollary}

\theoremstyle{definition}
\newtheorem{example}[thm]{Example}
\newtheorem{definition}[thm]{Definition}
\newtheorem{defi}[thm]{Definition}

\newtheorem{rem}[thm]{Remark}
\newtheorem*{rem*}{Remark}

\def\t{\times}

\newcommand{\Z}{\mathbb{Z}}

\newcommand{\N}{\mathbb{N}}

\newcommand{\Tor}{\text{Tor}}
\newcommand{\diff}{\partial}

\title[Kunneth Theorems for Vietoris-Rips Homology]{Kunneth Theorems for Vietoris-Rips Homology }
\author{Antonio Rieser}
\email{antonio.rieser@cimat.mx}
\author{Alejandra Trujillo-Negrete}
\address{Centro de Investigaci\'on en Matem\'aticas, A. C.  Jalisco S/N, Col. Valenciana CP: 36023 Guanajuato, Gto, México}
\email{alejandra.trujillo@cimat.mx}
\thanks{Research supported in part by Cat\'{e}dras CONACYT 1076, the US Office
of Naval Research Global, and the Southern Office of Aerospace Research and Development of the US Air Force Office of Scientific Research}

\begin{document}

\begin{abstract}
	We prove a Kunneth theorem for the Vietoris-Rips homology and
        cohomology of a semi-uniform space. We then interpret this result for
        graphs, where we show that the Kunneth theorem holds for graphs with
        respect to the strong graph product. We finish by computing the
        Vietoris-Rips cohomology of the torus endowed with diferent
        semi-uniform structures.
\end{abstract}

\maketitle
 
\section{Introduction}

The Vietoris-Rips complex was first defined by Vietoris in 1927 as a way to
obtain homology groups from metric spaces \cite{Vietoris_1927}, and, somewhat
later, it began to be used in the study of hyperbolic groups
\cite{Gromov_1987}. With the rise of topological data analysis in the
last fifteen years, the Vietoris-Rips homology has become a
computational, as well as theoretical, tool, and, indeed, it has become the standard invariant used 
in the homological analysis of data, in addition to its natural importance in the homological analysis of 
networks and graphs. Despite this surge of popularity, however, relatively little is known about the
properties of Vietoris-Rips homology, and, until recently, even many basic
results on the Vietoris-Rips homology had not been established. In a previous
article \cite{Rieser_2020}, the first author introduced a construction
of the Vietoris-Rips homology for semi-uniform spaces and proved a variant of
the Eilenberg-Steenrod axioms adapted to this context.

Studying the Vietoris-Rips homology and cohomology from the point of view of
semi-uniform spaces has many advantages over working directly with clique
complexes. The first is that setting up the Vietoris-Rips homology in this way enables one to use variations of
classical topological arguments in settings that are very far from
topological. Indeed, the results in this paper provide one example of this procedure. This, in turn, leads to
another major advantage of working in the category of semi-uniform spaces. As
the first author showed in \cite{Rieser_arXiv_VR_2020}, 
the Vietoris-Rips homology and cohomology are functors which 
act directly on semi-uniform spaces and uniformly continuous
maps between them. Since there exist non-trivial uniformly continuous maps between
many different kinds of spaces of interest, in
particular, between topological spaces, metric spaces with a preferred scale, and graphs,
this functoriality enables the construction of homomorphisms which can
be used to compare the resulting Vietoris-Rips homology
and cohomology groups across different classes of spaces. This is much more difficult to do if
one simply studies the Vietoris-Rips homology in each setting independently of the others.
Finally, using the constructions of semi-uniform spaces from \v Cech
closure spaces given in \cite{Rieser_arXiv_VR_2020}, we are able to see that the Vietoris cohomology of
Dowker \cite{Dowker_1952}, the metric cohomology of Hausmann \cite{Hausmann_1995}, 
the Vietoris-Rips cohomology on metric spaces with a preferred scale 
introduced in \cite{Rieser_arXiv_VR_2020}, and the cohomology of the flag complex 
of a graph are all instances of the same construction on different semi-uniform spaces.

In this article, we
continue the development of the algebraic topology of semi-uniform spaces
started in \cite{Rieser_arXiv_VR_2020}, first giving an alternate definition of
the Vietoris-Rips complex of a relation using simplicial
sets, and then studying the Vietoris-Rips homology
and cohomology of products of semi-uniform spaces with the goal of
establishing K\"unneth theorems in this context. For the
Vietoris-Rips homology as defined here, we will see that, while the Kunneth
theorem holds for Vietoris-Rips cohomology on arbitrary semi-uniform spaces, it is not true in general for
Vietoris-Rips homology. Nonetheless, we are able to show that it does hold for semi-uniform spaces induced
by graphs, which, in turn, implies a Kunneth theorem for the classical Vietoris-Rips
homology of graphs, which is the case of most interest to applications. Note that, while
it is well-known that the Kunneth theorem is false for the Vietoris-Rips homology using Cartesian products of 
graphs, by translating the problem into the setting of semi-uniform spaces, we
see that that one should use the 
\emph{strong graph product} instead. For a homology theory on graphs which
satisfies the Kunneth theorem with respect to 
the Cartesian graph product, see \cite{Grigoryan_et_al_2017}.

A Kunneth theorem for the classical Vietoris-Rips homology on
metric spaces with respect to the maximum metric on the product may 
also be deduced by applying the Kunneth theorem for simplicial
complexes to the isomorphism in Proposition 10.2 in
\cite{Adamaszek_Adams_2017}. However, with the
exception of the cases treated in \cite{Rieser_2020}, it 
remains unclear for which cases the the classical and
the semi-uniform Vietoris-Rips homology theories coincide. There has
also been some recent work on Kunneth theorems in persistent homology
\cites{Gakhar_Perea_arXiv_2019, Bubenik_Milicevic_arXiv_2019}, in which an expression for the persistent
homology of a product is obtained, given a filtered complex constructed from a
category whose homology has a Kunneth formula. The main
contribution of this article is that, by constructing the Vietoris-Rips homology
and cohomology
in the more general context of semi-uniform spaces, we are able to treat 
Kunneth theorems for the Vietoris-Rips cohomology of graphs, metric spaces, and
even topological spaces as
particular instances of the same theorem. Combined with the results of
\cite{Rieser_arXiv_VR_2020}, this allows for the computation of the Vietoris-Rips
cohomology of products of spaces which are semi-uniformly homotopy equivalent to spaces whose
Vietoris-Rips cohomologies are known (such as such as $(S^1 \times
S^1,d_{max})$, studied in \cite{Adamaszek_Adams_2017}), something which is unattainable using current techniques.

\section{Semi-uniform spaces and the Vietoris-Rips complex}

In this section, we recall the definition of semi-uniform spaces, which will be our main object of study. We begin with a few preliminary definitions.

\subsection{Semi-Uniform Spaces} 

\begin{defi}Let $U\subset X\times X$. We define $$U^{-1} \coloneqq\{(y,x)\mid
(x,y)\in U\}.$$ \end{defi}

\begin{definition}
	\label{def:Filter}
	Let $X$ be a set, and let $\mathcal{F}$ be a non-empty collection of subsets of $X$ with $\emptyset \notin \mathcal{F}$. We say that $\mathcal{F}$ is a \emph{filter} iff
	\begin{enumerate}
		\item $U \in \mathcal{F}$ and $U \subset V \implies V \in \mathcal{F}$, and
		\item \label{item:Filter condition 2} $U,V \in \mathcal{F} \implies U\cap V \in \mathcal{F}$.
	\end{enumerate}
\end{definition}

\begin{rem}
	The condition that $\mathcal{F}$ does not contain the empty set is
        occasionally additional in the literature, and such filters are
        sometimes called \emph{proper filters}. Since we will only be dealing
        with such filters, we see no need to make such a distinction here.
        Also, note that the requirement that $\emptyset \notin \mathcal{F}$
        combined with Condition \ref{item:Filter condition 2} of Definition
        \ref{def:Filter} above implies that the intersection of any finite
        collection of sets in $\mathcal{F}$ is nonempty.
\end{rem}

\begin{definition}
	\label{def:Semi-uniform}
     Let $X$ be a set. We say that a filter $\mathcal{U}$ on the product $X \times X$ is a \emph{semi-uniform structure on $X$} iff
     \begin{enumerate}
     	\item Each element of $\mathcal{U}$ contains the diagonal, i.e. $\Delta \subset U$ for all $U \in \mathcal{U}$
     	\item \label{item:Semi-uniform condition 2} If $U \in \mathcal{U}$, then $U^{-1}$ contains an element of $\mathcal{U}$
     \end{enumerate}
 	The pair $(X,\mathcal{U})$, consisting of a set $X$ and a semi-uniform structure $\mathcal{U}$ on $X$, is called a \emph{semi-uniform space}.	
\end{definition}

\begin{rem}
	Note that, since $\mathcal{U}$ is a filter, Condition \ref{item:Semi-uniform condition 2} in Definition \ref{def:Semi-uniform} is equivalent to the condition that $U^{-1} \in \mathcal{U}$.
\end{rem}

We give several important examples of semi-uniform spaces.

\begin{example}
    \label{ex:SU from graph}
Let $G = (V,E)$ be an undirected graph. We denote by $\sigma(G) \coloneqq
(V_G,\mathcal{E}_G)$ the semi-uniform space where $V_G = V$ and $\mathcal{E}_G
= [E\cup \Delta_{V_G}]$, i.e. $\mathcal{E}_G$ is the filter of subsets of $V_G
\times V_G$ generated by the set $E\cup \Delta_{V_G} \subset V_G \times V_G$,
where $\Delta_{V_G}$ is the diagonal in $V_G \times V_G$. 

Since $G$ is undirected in the above definition, $E = E^{-1}$, it follows from
Theorem 23.A.4 in \cite{Cech_1966}, that
$\sigma(G)$ is a semi-uniform space.
\end{example}
\begin{example}
    \label{ex:SU from metric}
Let $(X,d)$ be a metric space. For every $q > 0$, we define 
    \begin{align*}
        U_{q} &\coloneqq \{(x,y) \subset X \times X \mid d(x,y) <
        q\}\\
        U_{\leq q} &\coloneqq \{(x,y) \subset X \times X \mid d(x,y)
        \leq q\}.
        \end{align*}
    Now fix an $r>0$, and define $\mathcal{U}_r$ to be the semi-uniform
    structure generated by the sets $U_{r+\epsilon}$ for all $\epsilon>0$, and
    define $\mathcal{U}_{\leq r}$ to be the semi-uniform structure generated by
    the single relation $U_{\leq r}$. These are semi-uniform structures by Theorem 23.A.4
    in \cite{Cech_1966}.

\end{example}
\begin{example}[See also \cite{Rieser_arXiv_VR_2020}]
    \label{ex:SU from closure}
    Let $X$ be a set. We define a \emph{\v Cech closure operator}
    $c:\mathcal{P}(X) \to \mathcal{P}(C)$ to be a map such that
\begin{enumerate}
    \item $c(\emptyset) = \emptyset$
    \item $A \subset c(A)$ for all $A \subset X$
    \item $ c(A\cup B) = c(A) \cup c(B)$ for all $A,B \subset X$
\end{enumerate}
The pair $(X,c)$ will be called a \v Cech closure space. We define the interior
of a set $A \subset X$ to be
\begin{equation*}
    i(A) = X - c(X-A).
\end{equation*}
A set $U \subset X$ is said to be a \emph{neighborhood} of $A \subset X$ iff $A
\subset i(U)$, and a collection $\mathcal{U}=\{U_\alpha\mid \alpha \in A\}$ is
said to be an \emph{interior cover} of $X$ iff $X = \cup_\alpha i(U_\alpha)$.

For every interior cover $\mathcal{U}$ of a \v Cech closure space $(X,c)$, we
define two relations $V_\mathcal{U}$ and $R_\mathcal{U}$ by
\begin{align*}
    V_\mathcal{U} &\coloneqq \{(x,y) \mid \exists U \in \mathcal{U} : x,y \in
    U\}\\
        R_\mathcal{U} &\coloneqq \{(x,y) \mid \exists U \in \mathcal{U} : (x
            \in i(U) \textnormal{ and } y \in U) \textnormal{or} (y \in i(U)
        \textnormal{ and } x \in U\}.
\end{align*}
Let $\mathcal{I}_c$ denote the collection of interior covers of $(X,c)$, and
let $\mathcal{V}_c$ and $\mathcal{R}_c$ be the filters on $X \times X$
generated by the collections $\{V_\mathcal{U}\}_{\mathcal{U} \in
\mathcal{I}_c}$ and $\{R_{\mathcal{U}} \}_{\mathcal{U} \in \mathcal{I}_c}$,
respectively. As shown in \cite{Rieser_arXiv_VR_2020}, $\mathcal{V}_c$ and
$\mathcal{R}_c$ are semi-uniform structures on the set $X$. 

In particular, this construction produces a semi-uniform space for each topological space
$(X,c_\tau)$, where $c_\tau$ is a closure operator with $c^2_\tau = c_\tau$.
Note that in this case $\mathcal{V}_{c_\tau} = \mathcal{R}_{c_\tau}$, since both are generated by open covers.
\end{example}

The above examples show how semi-uniform spaces are a simultaneous
generalization of graphs, metric spaces with a preferred scale, and closure
spaces. We refer to \cites{Cech_1966, Rieser_arXiv_VR_2020} for further details on
semi-uniform spaces. 

\subsection{The Vietoris-Rips Homology and Cohomology of a Semi-Uniform Space}
\label{subsec:VR homology}
Let $(X,\mathcal{U})$ be a semi-uniform space and $U \in \mathcal{U}$.  We
define a simplicial set $X^U$ by setting $X^U_0
\coloneqq\{x\in X\} = X$, and, for $n \in \N$, we define 
\begin{equation*}
    X^U_n=\{(x_0,\dots,x_n)\mid (x_i,x_j)\in U) \; \forall i < j\},
\end{equation*}
with  functions 
\begin{align*} 
\diff_i&\colon  X^U_n\to X^U_{n-1}, \text{ where }
\diff_i(x_0,x_1,...,x_n)=(x_0,x_1,\dots,\hat{x}_i,\dots ,x_n), \\  
s_i&\colon X^U_n\to X^U_{n+1}, \text{ where }
s_i(x_0,x_1,...,x_n)=(x_0,x_1,\dots,x_i,x_i,\dots,x_n).
\end{align*}
We emphasize that the $x_i$ need not be distinct points of $X$, and, in
particular, since the diagonal $\Delta$ is a member of every $U \in
\mathcal{U}$, the $X^U_n$ will contain elements with $x_i = x_j$, $i\neq j$.

\begin{definition}
	We call the simplicial set $X^U$ the \emph{Vietoris-Rips complex of the pair $(X,U)$}.
\end{definition}

For each $n \in \Z$, let $C_n(X^U)$ be the graded free abelian group generated by the elements of the sets $X^U_n$, and let $C_*(X^U) = \oplus_{n \in \Z} C_n(X^U)$. Define a differential $d_n:C_n(X^U) \to C_{n-1}(X^U)$ by
\begin{equation*}
d_n = \sum_{i = 0}^n (-1)^i \diff_i
\end{equation*}
With these definitions, $(C_*(X^U),d)$ is now a chain complex, and we denote its homology  by $H_*(X^U)$. 

Let $\mathcal{G}$ be an $R$-module for a commutative ring $R$. Then $(C_*(X^U) \otimes_{\Z} \mathcal{G},d
\otimes \mathbf{1})$ is a chain complex over $R$, and we denote its homology, the
homology of $X^U$ with coefficients in $\mathcal{G}$, by
$H_*(X^U;\mathcal{G})$. Similarly, a $q$-dimensional cochain $f\in
C^q(X^U;\mathcal{G})$ is defined as a homomorphism $f\colon C_q(X^U)\otimes_\Z
R \to \mathcal{G}$ and the coboundary is given by 
\begin{equation*}
(\delta f) (\sigma)= \sum_{i=0}^{q+1}(-1)^i f(\diff_i \sigma)
\end{equation*}
for each $(q+1)$-simplex $\sigma $ of $X^U$. This leads to the cohomology groups
$H^*(X^U; \mathcal{G})$. 

We now consider maps between the simplicial sets corresponding to different elements of a semi-uniform structure $\mathcal{U}$.

\begin{proposition}
For $U,V \in \mathcal{U}$, $V \subset U$, $X^V$ is a sub-simplicial set of
$X^U$. Denote the inclusion map by $\phi_{VU}:X^V \hookrightarrow X^U$. For $W \subset V \subset U$, we have $\phi_{VU} \circ
\phi_{WV} = \phi_{WU}$.
\end{proposition}

\begin{proof}
Since $V \subset U$, if $\sigma = (x_0,\dots,x_n) \in X_n^V$, then, by
definition, $(x_i,x_j) \in V$ for all $i < j$. Therefore, $(x_i,x_j) \in U$ for
all $i < j$, and $\sigma \in X_n^U$. Define $\phi_{VU}(\sigma) = \sigma$. Since
this is both a simplicial map and an inclusion, the final statement follows, and the proof is complete.
\end{proof}

We define a partial order $\leq$ on the semi-uniform structure $\mathcal{U}$ by
writing $U \leq V$ iff $V \subseteq U$. Furthermore, since $\mathcal{U}$ is a
filter, for any $U,V \in \mathcal{U}$, $W = U \cap V \in \mathcal{U}$, and
therefore $\mathcal{U}$ with this partial order is a directed set. The induced
maps $\phi_{UV*}:H_*(X^V) \to H_*(X^U)$ and $\phi^*_{UV}:H^*(X^U) \to H^*(X^V)$
make $(H_*(X^U),\phi_{UV*},\mathcal{U})$ and $(H^*(X^U),
\phi_{UV}^*,\mathcal{U})$ into inverse and direct systems of abelian groups, respectively. We finally define the Vietoris-Rips homology and cohomology of a semi-uniform space $(X,\mathcal{U})$ to be
\begin{align*}
H^{VR}_*(X,\mathcal{U}) &= \varprojlim H_*(X^U) \\
H^*_{VR}(X,\mathcal{U}) &= \varinjlim H^*(X^U).
\end{align*}
We will typically suppress the semi-uniform structure $\mathcal{U}$ when it is unambiguous.

\begin{rem} Note that, since the ordered and unordered simplicial chain
    complexes are chain equivalent (by \cite{Spanier_1966}, Theorem 4.3.8), and since degenerate
    simplices in the simplicial set do not contribute to the resulting homology and
    cohomology (by \cite{MacLane_1975}, VIII.6), the Vietoris-Rips homology and cohomology defined here
for a semi-uniform space $(X,\mathcal{U})$ are isomorphic to the those defined in
\cite{Rieser_2020}.
\end{rem}

\begin{example}
    Let $G = (V,E)$ be a graph. The Vietoris-Rips homology and cohomology of
    $\sigma(G)=(V,\mathcal{E}_G)$ defined in Example \ref{ex:SU from graph} are the simplicial homology and
    cohomology, respectively, of the clique complex
    of the graph, by Theorem \ref{thm:VR graph} below.
\end{example}

\begin{example}
    Let $(X,d)$ be a metric space, and let $(X,\mathcal{U}_{r})$ be the
    semi-uniform spaces constructed in Example \ref{ex:SU from metric}. Since
    the sets $U_\epsilon$, $\epsilon > 0$ are cofinal in $\mathcal{U}_0$, it
    follows that the Vietoris-Rips cohomology of $(X,\mathcal{U}_0)$ is isomorphic to the \emph{metric
    cohomology} studied by Hausmann in \cite{Hausmann_1995}.

\end{example}

\begin{example}
Let $(X,c_{\tau})$ be a topological closure space (where $c_\tau^2 = c_\tau$),
and let $(X,\mathcal{V}_{c_\tau})$ be the corresponding semi-uniform space from
Example \ref{ex:SU from closure}. By results in Dowker \cite{Dowker_1952},
the Vietoris-Rips homology and cohomology of $(X,\mathcal{V}_{c_\tau})$ are isomorphic to the \v Cech
homology and cohomology of the topological space $(X,\tau)$, where the topology
$\tau$ is generated by
$c_\tau$, i.e. a set $U \subset X$ is open in $(X,\tau)$ iff
$c_\tau^2(X - U) =c_\tau(X-U)$.
\end{example}

\section{Products}

In this section, we recall the definitions of products for semi-uniform spaces and simplicial sets, and we then prove a theorem relating the products of Vietoris-Rips complexes which will be the basis for the Kunneth Theorems to follow.

\subsection{Products of Semi-Uniform Spaces}

\begin{definition}
Let $\{(X_a,\mathcal{U}_a)\}_{a \in A}$ be a family of semi-uniform spaces indexed by $A$. Let
$X \coloneqq \Pi_{a \in A} X_a$ be the Cartesian product of the sets, and we
denote by $\mathcal{U}$ the filter generated by the subsets of $X \times X$ of the form
\begin{equation}
\label{eq:Product filter base}
\{ (x,y) \mid (x,y) \in X \times X, a \in F \implies (\pi_a x, \pi_a y) \in U_a\},
\end{equation}
where $F \subset A$ is some finite subset of $A$, $U_a \in \mathcal{U}_a$ for
each $a \in F$, and $\pi_a \colon X \to X_a$ is the projection to the $a$-th
coordinate. 
\end{definition}

\begin{prop}
$(X,\mathcal{U})$ defined as above is a semi-uniform space.
\end{prop}

\begin{proof}
First, let $(x,x) \in \Delta_{X \times X}$. Since the $\Delta_{U_a} \subset
U_a$ for every $U_a \in \mathcal{U}_a$ and every $a \in A$, it follows that
$(x,x)$ is in every element of $\mathcal{U}$. Since $(x,x) \in \Delta_{X \times X}$
was arbitrary, $\Delta_{X \times X}$ is contained in every element $\in \mathcal{U}$.

Suppose now that $U \in \mathcal{U}$. Then $U$ contains a set $V$ of form
(\ref{eq:Product filter base}) above. Since, for any $U_a \in \mathcal{U}_a$ we
have that
$U_a^{-1} \in \mathcal{U}_a$, we see from $V_a^{-1} \subset U_a^{-1}$
that $V^{-1}$ is of the form (\ref{eq:Product filter base}) as well. However, $V^{-1} \subset U^{-1}$, so $U^{-1} \in \mathcal{U}$, and the proof is complete.
\end{proof}

\begin{definition} We call $(X,\mathcal{U})$ the \emph{semi-uniform product of }$\{(X_a,\mathcal{U}_a)\}_{a \in A}$. We will sometimes denote $(X,\mathcal{U})$ as $(\Pi_{a \in A} X_a,\Pi_{a \in A} \mathcal{U}_a)$.
\end{definition}

\subsection{Products of Simplicial Sets}

\begin{definition}
	Let $\Sigma$ and $\Sigma'$ be simplicial sets. Then the product $\Sigma \times \Sigma'$ is given by
	\begin{equation*}
            (\Sigma \times \Sigma')_n \coloneqq \{ (\sigma,\sigma') \mid \sigma \in \Sigma_n, \sigma' \in \Sigma'_n\}
	\end{equation*}
\end{definition}
For any simplicial sets $\Sigma$ and $\Sigma'$, the Eilenberg-Zilber theorem
\cite{Eilenberg_Zilber_1953} gives a quasi-isomorphism between the chain
complexes $C_*(\Sigma \times \Sigma')$ and $C_*(\Sigma) \otimes C_*(\Sigma')$.
In order to establish the Kunneth theorems for Vietoris-Rips homology, we must
further establish a relationship between the Vietoris-Rips complexes of
products of semi-uniform spaces on the one hand, and the products of
Vietoris-Rips complexes of semi-uniform spaces on the other. This follows
easily from the respective definitions, and is accomplished in the following
theorem. First, however, we make the following remark.

\begin{rem*}
Let $(X,\mathcal{U})$ and $(Y,\mathcal{V})$ be semi-uniform spaces and let  $U \in \mathcal{U}$ and $V \in \mathcal{V}$.   Observe that the set $U\t V\subset (X\t X) \t (Y\t Y) $ can be seen as  a subset of   $(X\t Y) \t (X\t Y)$ via the isomorphism 
\begin{align*}
\psi\colon  (X\t Y)\t(X\t Y)&\to  (X\t X) \t (Y\t Y) \\ ((x_1,y_1),(x_2,y_2))&\mapsto ((x_1,x_2), (y_1,y_2)
\end{align*} 
Note that the expression $(X \times Y)^{U\times V}$ is an abuse of notation, and is, more precisely,
$(X \times Y)^{\psi^{-1}(U\times V)}$. Since $U \times V \cong \psi^{-1}(U \times V)$, however,
we will use the first notation in place of the second throughout.
\end{rem*}
	\begin{theorem}
		\label{thm:Simplicial set of product}
	Let $(X,\mathcal{U})$ and $(Y,\mathcal{V})$ be semi-uniform spaces. Then, for any $U \in \mathcal{U}$ and $V \in \mathcal{V}$, we have
	\begin{equation*}
	(X \times Y)^{U \times V} \cong X^U \times Y^V
	\end{equation*}
	\end{theorem}

	\begin{proof}
		Let $\phi\colon (X \times Y)^{U\times V}_k \to (X^U \times Y^V)_k $ be the map given by
		\begin{equation*}
		\phi ((x_0,y_0),\dots, (x_k,y_k)) \coloneqq ((x_0,\dots,x_k),(y_0,\dots,y_k)).	
		\end{equation*}
		By definition, $((x_0,y_0),\dots, (x_k,y_k)) \in (X \times Y)^{U\times V}$ iff, for all $i,j \in \{1,\dots,k\}$, one of the following holds: 
		\begin{enumerate}
		\item $\psi((x_i,y_j),(x_j,y_j)) \in U \times V$
		\item $x_i = x_j$ and $(y_i,y_j) \in V$
		\item $(x_i,x_j) \in U$ and $y_i = y_j$
		\end{enumerate}
		This, in turn, is true iff $(x_0,\dots,x_k) \in X^U$ and $(y_0,\dots,y_k) \in Y^V$. Therefore, $\phi$ gives an isomorphism between $(X \times Y)^{U \times V}$ and $X^U \times Y^V$, and the proof is complete.
	\end{proof}

\section{The Kunneth Theorems}
\label{sec:Kunneth theorems semi-uniform}

The Kunneth theorems for the Vietoris-Rips homology now follow from the above
product relations, the Kunneth theorems for simplicial sets, and the properties
of exact sequences in direct and inverse limits. 
We begin with the
results for Vietoris-Rips cohomology, where we have a Kunneth formula in
general. 

In order to establish this, we will first require the following lemma. We begin
by recalling the following definition.

\begin{definition}
We say that a graded module $\{C_q\}$ is said to be of \emph{finite type} if it is finitely generated for
every $q$. 
\end{definition}

\begin{rem} We will sometimes write the torsion product
    $\Tor_1(A,B)$ as $A*B$.
\end{rem}
\begin{lem}
	\label{lem:One-set Kunneth for Cohomology} 
	Let $(X,\mathcal{U})$ and $(Y,\mathcal{V})$ be semi-uniform spaces, $U
        \in \mathcal{U}$, and $V \in \mathcal{V}$. Let $\mathcal{G}$ and
        $\mathcal{G}'$ be modules over a principal ideal domain $R$ such that
        $\mathcal{G}*\mathcal{G}'= 0$. If $H_*(X^U;R)$ and $H_*(Y^V;R)$ are of finite type or if
        $H_*(Y^V;R)$ is of finite type and $\mathcal{G}'$ is finitely generated, then for
        any $q \in \Z$ there is a natural short-exact sequence
	\begin{align*}
            0 \to \oplus_{i+j=q}& H^i_{VR}(X^U;\mathcal{G})\otimes
            H^j_{VR}(Y^V;\mathcal{G}') 
            \to H^q_{VR}((X \times Y)^{U \times V};\mathcal{G} \otimes
            \mathcal{G}')\to\\
                                &\to \oplus_{i+j=q+1}
                                \textnormal{Tor}_1(H^i_{VR}(X^U;\mathcal{G}),H^j_{VR}(Y^V;\mathcal{G}'))
        \to 0.
	\end{align*}
        Furthermore, this sequence splits, but not canonically.
\end{lem}
\begin{proof}
	Theorem \ref{thm:Simplicial set of product} gives $C_*(X\times Y,U
        \times V) \cong C_*(X^U \times Y^V)$, and therefore $C^*(X\times Y,U
        \times V) \cong C^*(X^U \times Y^V)$, and from the Eilenberg-Zilber
        theorem, we have $C^*(X^U \times Y^V) \simeq C^*(X^U)\otimes C^*(Y^V)$,
        where $\simeq$ indicates cochain-homotopy equivalence. The result now
        follows from the Kunneth theorem for cochain complexes from
        \cite{Spanier_1966}, Theorem 5.5.11.
\end{proof}

We now give the Kunneth theorem for Vietoris-Rips cohomology.

\begin{theorem}
	\label{thm:Kunneth theorem VR cohomology}
	Let $(X,\mathcal{U})$ and $(Y,\mathcal{V})$ be semi-uniform spaces, and
        let $\mathcal{G}$ and $\mathcal{G}'$ be modules over a principal ideal domain $R$ such that
        $\mathcal{G}*\mathcal{G}'= 0$. Suppose that $\mathcal{U}$ and $\mathcal{V}$ are generated by
        collections $\{U_\alpha\}_{\alpha \in A}$ and $\{V_\beta\}_{\beta \in
        B}$ such that $H_*(X^{U_\alpha};R)$ and $H_*(Y^{V_\beta};R)$ are of
        finite type for all $\alpha \in A$ and $\beta \in B$, or such that
        $H_*(Y^{V_\beta};R)$ is
        of finite type for $\beta \in B$ and $\mathcal{G}'$ is finitely generated. 
        Then for any $q \in \Z$ there is a natural short-exact sequence
        \begin{align*}
            0 \to \oplus_{i+j=q}& H^i_{VR}(X;\mathcal{G})\otimes
            H^j_{VR}(Y;\mathcal{G}') 
            \to H^q_{VR}(X \times Y;\mathcal{G} \otimes \mathcal{G}')\to\\
                                &\to \oplus_{i+j=q+1}
                                \textnormal{Tor}_1(H^i_{VR}(X;\mathcal{G}),H^j_{VR}(Y;\mathcal{G}')) \to 0.
	\end{align*}
        Furthermore, this sequence splits, but not canonically.
\end{theorem}
\begin{proof}
	We first note that the family of sets in $\mathcal{U} \times
        \mathcal{V}$ of the form $U \times V$ are cofinal in $\mathcal{U}
        \times \mathcal{V}$. The theorem now follows from Proposition
        \ref{lem:One-set Kunneth for Cohomology}, the fact that $Tor_1$ and $\otimes$ commute with direct limits, and that the direct limit is an exact functor.
\end{proof}

For homology, the Kunneth theorem does not hold in general, due to the failure
of exactness of the inverse limits of exact sequences. Nonetheless, in some
special cases, the exact sequences related to elements $U \times V \in
\mathcal{U}\times \mathcal{V}$ prove to be useful. We give two such situations
below, beginning, as above with the following lemma.

\begin{lem} 
	\label{thm:One-set Kunneth homology} 
	Let $(X,\mathcal{U})$ and $(Y,\mathcal{V})$ be semi-uniform spaces, and
        let $\mathcal{G}$ and $\mathcal{G}'$ be modules or a principal ideal domain such that
        $\mathcal{G}*\mathcal{G}'= 0$. Then for any $U \in \mathcal{U}$, $V \in \mathcal{V}$, and
        $q \in \Z$, there is a natural short-exact sequence
	\begin{align*}
            0 \to \oplus_{i+j=q}& H_i^{VR}(X^U;\mathcal{G})\otimes
            H_j^{VR}(Y^V;\mathcal{G}') 
            \to H_q^{VR}((X \times Y)^{U \times V};\mathcal{G} \otimes
            \mathcal{G}')\to\\
                                &\to \oplus_{i+j=q-1}
                                \textnormal{Tor}_1(H_i^{VR}(X^U;\mathcal{G}),H_j^{VR}(Y^V;\mathcal{G}')) \to 0.
	\end{align*}
        Furthermore, this sequence splits, but not canonically.
\end{lem}
\begin{proof}
	Theorem \ref{thm:Simplicial set of product} gives $C_*(X\times Y,U
        \times V) \cong C_*(X^U \times Y^V)$, and from the Eilenberg-Zilber
        theorem \cite{Eilenberg_Zilber_1953}, we have $C_*(X^U \times Y^V)
        \simeq C_*(X^U)\otimes C_*(Y^V)$, where $\simeq$ indicates
        chain-homotopy equivalence. The result now follows from the Kunneth
        theorem for chain complexes (Theorem 5.3.4 in \cite{Spanier_1966}).
\end{proof}
\begin{theorem}
	\label{thm:Kunneth homology - max}
	 Let $(X,\mathcal{U})$ and $(Y,\mathcal{V})$ be semi-uniform spaces,
         and let $\mathcal{G}$ and $\mathcal{G}'$ be modules over a principal
         ideal domain such that $\mathcal{G} * \mathcal{G}'=0$. Suppose $U^*$
         and $V^*$ are maximal in $\mathcal{U}$ and $\mathcal{V}$,
         respectively, ordered by inclusion. Then $U^* \times V^*$ is maximal
         in $\mathcal{U} \times \mathcal{V}$, and for all $q \in \Z$, there
         exists a natural short-exact sequence
	\begin{align*}
            0 \to \oplus_{i+j=q}& H_i^{VR}(X;\mathcal{G})\otimes
            H_j^{VR}(Y;\mathcal{G}') 
            \to H_q^{VR}(X \times Y;\mathcal{G}\otimes\mathcal{G}')\to\\
                                &\to \oplus_{i+j=q-1}
                                \textnormal{Tor}_1(H_i^{VR}(X;\mathcal{G}),H_j^{VR}(Y;\mathcal{G}'))
                                \to 0.
	\end{align*}
        Furthermore, this sequence splits, but not canonically.
\end{theorem}

\begin{proof}
	We first show that $U^* \times V^*$ is maximal in $\mathcal{U}\times \mathcal{V}$. Suppose that there exists a $W \in \mathcal{U} \times \mathcal{V}$ with $W \subsetneq U^* \times V^*$. Then there exist $U \in \mathcal{U}$ and $V \in \mathcal{V}$ such that $U \times V \subset W \subsetneq U^* \times V^*$, a contradiction. Therefore $U^* \times V^*$ is maximal in $\mathcal{U} \times \mathcal{V}$.
	For any semi-uniform space $(X,\mathcal{U})$, we have, by definition,
        $H_*^{VR}(X,\mathcal{U};\mathcal{G}) = \varprojlim
        H_*^{VR}(X^U;\mathcal{G})$. If $U^*$ is maximal in $\mathcal{U}$, then
        $H_*(X^{U^*};\mathcal{G})$ is cofinal in the inverse system
        $\{H_{*}(X^U;\mathcal{G}),\phi_{UV*},\mathcal{U}\}$, from which it follows that
	\begin{equation*}
            H_*^{VR}(X,\mathcal{U};\mathcal{G}) = \varprojlim
            H_*^{VR}(X^U;\mathcal{G}) = H_*(X^{U^*};\mathcal{G})
	\end{equation*}
	as desired. 
	Putting these together, the result now follows from Theorem \ref{thm:One-set Kunneth homology}.
\end{proof}

Although the short exact sequence doesn't hold in general for Vietoris-Rips
homology, if the torsion term vanishes on a cofinal subset of the bases for the 
product semi-uniform structure, we may still conclude that there is an isomorphism 
of the respective homology groups, as we see from the following.

\begin{theorem}
    \label{thm:Homology iso}
    Let $(X,\mathcal{U})$ and $(Y,\mathcal{V})$ be semi-uniform spaces, and let
    $\mathcal{G}$ and $\mathcal{G}'$ be modules over a principal ideal domain
    such that $\mathcal{G}*\mathcal{G}'= 0$.
    Suppose that $\mathcal{F}'\subset \mathcal{U}\times \mathcal{V}$ is a cofinal collection of sets in
    of the form $U \times V$, $U \in \mathcal{U}$, $V \in \mathcal{V}$ such
    that, for any $U \times V \in \mathcal{F}$, $\oplus_{i+j=q-1}
    \text{Tor}_1(H_i^{VR}(X;\mathcal{G}),H_j^{VR}(Y;\mathcal{G}')) = 0$. Then
    \[\oplus_{i+j=q} H_i^{VR}(X;\mathcal{G})\otimes H_j^{VR}(Y;\mathcal{G}') 
    \cong H_q^{VR}(X \times Y;\mathcal{G} \otimes \mathcal{G}')
\in \mathcal{V}.\]
\end{theorem}

\section{Kunneth Theorems for Vietoris-Rips homology on graphs}
\label{sec:Graphs}
	
In this section, we apply the general Kunneth formulae on semi-uniform spaces from Section
\ref{sec:Kunneth theorems semi-uniform} to prove the Kunneth Theorem for the Vietoris-Rips homology on graphs. We begin with the following construction. All graphs are undirected.

\begin{definition}
    Let $G = (V,E)$ and $G' = (V',E')$ (undirected) graphs. Then the strong
    graph product $G \boxtimes G' = (V \boxtimes V', E \boxtimes E')$ is defined by
	\begin{align*}
		(v,v')& \in V \boxtimes V' \iff v \in V, v'\in V'\\
		((v_0,v'_0), (v_1,v'_1)) &\in E \boxtimes E' \iff \text{One of the following holds:}\\
		& (1)\;\; (v_0,v_1) \in E \text{ and } (v'_0,v'_1) \in E'\\
		& (2)\;\; v_0 = v_1 \text{ and } (v'_0,v'_1) \in E'\\
		& (3)\;\; (v_0,v_1) \in E \text { and } v'_0=v'_1 
	\end{align*}
\end{definition}

\begin{theorem}
    \label{thm:VR graph}
    Let $H_*^{VR}(G;\mathcal{A})$ denote the Vietoris-Rips homology of the
    undirected graph $G = (V,E)$ with coefficients in an abelian group
    $\mathcal{A}$. Then $H_*^{VR}(G;\mathcal{A}) \cong
    H_*^{VR}(\sigma(G);\mathcal{A})$, where $\sigma(G) =  (V, \mathcal{E}_G)$ is the semi-uniform
    space defined in Example \ref{ex:SU from graph}.
\end{theorem}
\begin{proof}
	By construction, $\mathcal{E}_G$ has a maximal element $E$ consisting
        of the edges of the graph $G$ and the diagonal $\Delta_V$, from which
        it follows that $H_*^{VR}(\sigma(G);\mathcal{A}) \cong
        H_*(V_G^E;\mathcal{A})$. However, $H_*(V_G^E;\mathcal{A})$ is exactly the homology of the simplicial set generated by the clique complex $\Sigma_G$ of $G$. Therefore,
	\begin{equation*}
            H_*^{VR}(G;\mathcal{A}) = H_*(\Sigma_G;\mathcal{A}) \cong
            H_*(V_G^E;\mathcal{A}) \cong H_*^{VR}(\sigma(G);\mathcal{A}),
	\end{equation*}
	and the proof is complete.
\end{proof}

In order to prove the Kunneth Theorem for undirected graphs, we will need the following

\begin{proposition}
	\label{prop:Product graph semi-uniform space}
	Let $G = (V,E)$ and $G' = (V',E')$ be graphs. Then 
        \[(V \times V', \mathcal{E}_{G \boxtimes G'}) =
            \sigma(G \boxtimes G') = \sigma(G) \times \sigma(G') =
        (V,\mathcal{E}_G) \times (V',\mathcal{E}_{G'}),\]
        where the product on the right is the
        product of semi-uniform spaces.
\end{proposition}
\begin{proof}
	We must show that $\mathcal{E}_G \times \mathcal{E}_{G'} = \mathcal{E}_{G \boxtimes G'}$. First, let 
        $U \in \mathcal{E}_G \times \mathcal{E}_{G'}$. Then, by construction, there exist sets $e \in \mathcal{E}_G$ and $e'\in 
        \mathcal{E}_{G'}$ such that $e \times e' \subset U$. By definition of $\mathcal{E}_G$ and $\mathcal{E}_{G'}$, however, 
        $E \subset e$ and $E' \subset e'$, so $E \times E' \subset e \times e'
        \subset U$. Since $E \times E'\in  \mathcal{E}_{G \boxtimes G'}$ and $
        \mathcal{E}_{G \boxtimes G'}$ is a filter, therefore 
        $U \in \mathcal{E}_{G \boxtimes G'}$, and we see that 
        $\mathcal{E}_{G} \times \mathcal{E}_{G'} \subset \mathcal{E}_{G \boxtimes G'}$. 

        Now suppose $U \in \mathcal{E}_{G \boxtimes G'}$. Then $E \times E'
        \subset U$, and therefore $U \in \mathcal{E}_G \times \mathcal{E}_{G'}$
        as well, by definition of the product semi-uniform structure. Therefore
        $\mathcal{E}_{G \boxtimes G'} \subset  \mathcal{E}_{G} \times
        \mathcal{E}_{G'}$. It follows that $\mathcal{E}_G \times \mathcal{E}_{G'} = \mathcal{E}_{G \boxtimes G'}$, and therefore $\sigma(G \boxtimes G') = \sigma(G) \times \sigma(G)$.
\end{proof}

\begin{theorem} Let $G=(V,E)$ and $G'=(V',E')$ be graphs, and let $\mathcal{M}$
    and $\mathcal{M}'$ be modules over a principal ideal domain such that
    $\mathcal{M}*\mathcal{M}'=0$. For every $q \in \Z$, there exist natural short exact sequences
	\begin{align*}
            0 \to \oplus_{i+j=q} H_i^{VR}(G;\mathcal{M})\otimes
            H_j^{VR}&(G';\mathcal{M}') 
            \to H_q^{VR}(G \boxtimes G';\mathcal{M}\otimes\mathcal{M}')\to\\
                    &\to \oplus_{i+j=q-1}
                    \textnormal{Tor}_1(H_i^{VR}(G;\mathcal{M}),H_j^{VR}(G;\mathcal{M}')) \to 0
	\end{align*}
\end{theorem}
\begin{proof}
	The theorem follows immediately from Proposition \ref{prop:Product
        graph semi-uniform space}, Theorem \ref{thm:Kunneth homology - max}, and the definition of the classical Vietoris-Rips complex for graphs.
\end{proof}

\section{Applications to Metric Spaces}
\label{sec:Metric spaces}

As an illustration of the above results, we examine the Vietoris-Rips
homology and cohomology of the torus with different semi-uniform structures. We
first recall the definition of the classical Vietoris-Rips complex on a metric
space and the corresponding relations which can be used to generate
semi-uniform structures.    

\begin{definition} Suppose that $(X,d)$ is a metric space and $r > 0$ is a real number. The
    Vietoris-Rips complex $VR_{<}(X;r)$ is the simplicial complex with
    vertex set $X$, where a finite subset $\sigma\subseteq X$ is a simplex if
    only if the $diam(\sigma) < r$. The Vietoris-Rips complex $VR_{\leq}(X;r)$ 
    is the simplicial complex with vertex set $X$, where a finite subset $\sigma\subseteq X$ is a simplex if
    only if the $diam(\sigma) \leq r$.   
\end{definition}

\begin{rem}
    \label{rem:Geometric realization}
    For a given $r>0$, note that $VR_{<}(X;r)$ is the geometric realization of
    the simplicial set $X^{U_r}$, and $VR_{\leq}(X;r)$ is the geometric
    realization of the simplicial set $X^{U_{\leq
    r}}$. They therefore have the same homology and cohomology groups.
\end{rem}
We now recall following theorem from \cite{Adamaszek_Adams_2017}. 

\begin{thm}[\cite{Adamaszek_Adams_2017}, Theorem 7.4]\label{thm:Circle}
    Denote the circle with unit circumference by $\mathbb{S}^1$, and consider
    $\mathbb{S}^1$ as a metric space with the geodesic distance. For $0 < r < \frac{1}{2} $, suppose that
    $\frac{1}{2l+1}<r < \frac{l+1}{2l+3}$ for some $l \in \{0,1,...\} $. Then we have a homotopy equivalence 
\begin{equation}
VR_{< }(\mathbb{S}^1;r)\simeq S^{2l+1}. \end{equation}
\end{thm}

The following computation now follows from Theorem \ref{thm:Circle}
and Theorem \ref{thm:Kunneth theorem VR cohomology}.

\begin{proposition}
    \label{prop:Torus}
    Let $0< r, r'\leq \frac{1}{2}$ with $\frac{1}{2l+1}<r < \frac{l+1}{2l+3}$ and
    $\frac{1}{2l'+1}<r' < \frac{l'+1}{2l'+3}$ for some $l,l' \in
    \{0,1,2,\dots\}$.
    Let $\mathbb{T}^2_{r,r'}$ denote the product semi-uniform space
    $\mathbb{T}^2_{r,r'}\coloneqq (\mathbb{S}^1,\mathcal{U}_r) \times
    (\mathbb{S}^1,\mathcal{U}_{r'})$.

    If $l \neq l'$, then

    \begin{equation*}
        H^q_{VR}(\mathbb{T}^2_{r,r'}) = 
        \begin{cases}
            \mathbb{Z}              & q= 0,2l+1,2l'+1,\text{ or }2(l+l'+1)\\
            \{ 0 \}                 & \text{otherwise.}
        \end{cases}
    \end{equation*}

    If $l = l'$, then
    \begin{equation*}
        H^q_{VR}(\mathbb{T}^2_{r,r}) = 
        \begin{cases}
            \mathbb{Z}                              & q= 0 \text{ or }2(2l+1)\\
            \mathbb{Z} \times \mathbb{Z}            & q = 2l+1\\
            \{ 0 \}                                 & \text{otherwise.}
        \end{cases}
    \end{equation*}
 \end{proposition}
 \begin{proof}
By Theorem \ref{thm:Kunneth theorem VR cohomology}, we have the following short exact sequence  
\begin{align*}
	0 \to \oplus_{i+j=q}& H^i_{VR}(\mathbb{S}^1,\mathcal{U}_r)\otimes
        H^j_{VR}(\mathbb{S}^1,\mathcal{U}_{r'}) 
        \to H^q_{VR}(\mathbb{T}^2_{r,r'})\to\\
	&\to \oplus_{i+j=q-1} \text{Tor}_1(H^i_{VR}(\mathbb{S}^1,
        \mathcal{U}_r),H^j_{VR}(\mathbb{S}^1, \mathcal{U}_{r'})) \to 0
	\end{align*}
        By Remark \ref{rem:Geometric realization} and Theorem \ref{thm:Circle} we have that
        \[H^i_{VR}(X^{U_{r+\epsilon}})\cong
H^i(VR_{<}(X;r+\epsilon)) \cong H^i(S^{2l+1})\] for all $\epsilon>0$ sufficiently
small. Since the sets $U_{r+\epsilon}$ are cofinal in $\mathcal{U}_r$, 
we obtain the exact sequence 
\begin{align*}
0 \to \oplus_{i+j=q}& H_i(S^{2l+1})\otimes H_j(S^{2l'+1}) 
\to H_q^{VR}(\mathbb{T}^2_{r,r'})\to\\
	&\to \oplus_{i+j=q-1} \text{Tor}_1(H_i(S^{2l+1}),H_j(S^{2l'+1})) \to 0.
	\end{align*}
Since the torsion term in the above exact sequence is trivial, the result
follows. 
\end{proof}

\begin{rem} Note that a similar result is true for the Vietoris-Rips homology
    of the torus by Theorem \ref{thm:Homology iso}.
\end{rem}

Now, suppose that $(X,d_X)$ and $(Y,d_Y)$ metric spaces and $(X\times Y, d)$ with
$d((x_1,y_1),(x_2,y_2))=max\{d_X(x_1,x_2), d_Y(y_1,y_2)\}$.  Then we have 
\begin{align*}
U_q^{X\times Y}&=
\{(z,w)\in (X\times Y)\times (X\times Y) \mid d(z,w) < q \}
\\
& =\{((x_1,y_1),(x_2,y_2)) \mid d((x_1,y_1),(x_2,y_2)) < q \} 
\\
&=\{((x_1,y_1),(x_2,y_2)) \mid max\{d_X(x_1,x_2), d_Y(y_1,y_2)\} < q \}
\\
&= \{((x_1,y_1),(x_2,y_2)) \mid d_X(x_1,x_2) < r \text{ and }d_Y(y_1,y_2) < q\}\\
&\cong U_q^X \times U_q^Y
\end{align*}
Note, too, that the $U^{X \times Y}_{r+\epsilon}= $ are cofinal in
$\mathcal{U}_r^{X\times Y}$. The following corollary follows immediately from
the above comments, Theorem \ref{thm:Kunneth theorem VR cohomology}, and
Proposition \ref{prop:Torus}.
\begin{cor}Let $(X,d_X)$ and $(Y,d_Y)$ be metric spaces, and let $d$ be the
    maximum metric on $X \times Y$. Let $\mathcal{G}$ and $\mathcal{G}'$ be
    modules of a principal ideal domain such that $\mathcal{G}*\mathcal{G}'=0$. Then we have the short exact sequence 
    \begin{align*}
        0 \to \oplus_{i+j=q}& H^i_{VR}(X,\mathcal{U}_r;\mathcal{G})\otimes
        H^j_{VR}(Y,\mathcal{U}_{r};\mathcal{G}') 
        \to H^q_{VR}(X \times Y, \mathcal{U}_r;\mathcal{G} \otimes \mathcal{G}') \to\\
	&\to \oplus_{i+j=q-1} \textnormal{Tor}_1(H^i_{VR}(X,
        \mathcal{U}_r;\mathcal{G}),H^j_{VR}(Y, \mathcal{U}_{r};\mathcal{G}')) \to 0.
    \end{align*}
    Furthermore, this sequence splits, but not canonically.
\end{cor}

Applying this to the case of the torus, we recover the following result, which
also follows from Proposition 10.2 in \cite{Adamaszek_Adams_2017}.

\begin{cor}Let $\mathbb{T}^2=\mathbb{S}^1\times \mathbb{S}^1$ with the maximum metric. 
Let  $0 < r\leq \frac{1}{2}$ with  $\frac{1}{2l+1}<r < \frac{l+1}{2l+3}$ for some $l=0,1,... $. Then 
\begin{align*}
    H_q^{VR}(\mathbb{T}^2,\mathcal{U}_r)&\cong  \oplus_{i+j=q}
    H^i_{VR}(\mathbb{S}^1,\mathcal{U}_r)\otimes
    H^j_{VR}(\mathbb{S}^1,\mathcal{U}_{r}) \\
                                        &\cong
\begin{cases} \mathbb{Z} & \text{ if } q\in\{0, 2(2l+1)\}\\
    \mathbb{Z} \times \mathbb{Z} & \text{ if } q = 2l+1\\ 
    0 & otherwise.
\end{cases}
\end{align*}  
\end{cor}

\bibliography{/home/antonio/Bib/all.bib}

\end{document}